\documentclass[12pt,bezier]{article}
\usepackage{amssymb,amsfonts,amsthm}
\usepackage{eucal,latexsym}
\newcommand{\be}[1]{\begin{equation}\label{#1}}
\newcommand{\ee}{\end{equation}}
\newcommand{\rav}{\langle\,\Sigma\mid\rr\,\rangle}
\newcommand{\pres}[2]{\langle\,#1\mid #2\,\rangle}

\newcommand{\iv}{^{-1}}

\newcommand{\bs}{w}

\newcommand{\topp}{\mathop{\mbox{\bf top}}}
\newcommand{\bott}{\mathop{\mbox{\bf bot}}}

\newcommand{\rr}{{\cal R} }
\newcommand{\pp}{{\cal P} }
\newcommand{\qq}{{\cal Q} }

\newcommand{\zz}{{\bf Z}}

\newcommand{\dg}{{\cal D}(\pp,\bs)}
\renewcommand{\wr}{\mathrel{\rm wr}}

\newcommand{\comp}{\mathop{\bf comp\,}}

\newcommand{\la}{\langle\,}
\newcommand{\ra}{\,\rangle}

\begin{document}

\title
{Rigidity properties of diagram groups}
\author{V.\,S.~Guba \thanks{The research of the first author is supported in part by the RFFI grant
96--01--00974.},\ \,M.\,V.\,Sapir \thanks{The research of the
second author is supported in part by the NSF grant DMS 9623284}}
\date{}

\maketitle

\begin{abstract} In this paper we establish a rigid connection
between two classical objects: the R.Thompson group (the group of all
piece-wise linear homeomorphisms of the unit interval with finitely many dyadic break points and all slopes powers of 2) and the Dunce hat (the topological space obtained from the triangle ABC by gluing AB, BC and AC).
We prove that a diagram group of a directed 2-complex contains a copy of
the R.Thompson group if and only if the 2-complex contains a copy of
the Dunce hut.
\end{abstract}

\setlength{\unitlength}{0.001in}
\newtheorem{thm}{\quad Theorem}
\newtheorem{lm}[thm]{\quad Lemma}
\newtheorem{cy}[thm]{\quad Corollary}
\newtheorem{df}[thm]{\quad Definition}
\newtheorem{ex}[thm]{\quad Example}
\newtheorem{rk}[thm]{\quad Remark}
\newtheorem{prob}{\quad Problem}

\section{Introduction}

The class of diagram groups was introduced by Meakin and Sapir in
1993. Kilibarda obtained the first results about diagram groups in
\cite{KilDiss}, \cite{Kil}. The theory was further developed in
\cite{GuSa97}, \cite{GuSa99}. It turned out that many important
groups (including the R.\,Thompson group $F$) are diagram groups.
On the other hand, diagram groups satisfy some interesting
properties, and there exists a deep similarity between
combinatorics on diagrams and combinatorics on words. Recent results
by D. Farley \cite{Farley} show that diagram groups act by isometries
on CAT(0)-spaces. This allowed him to prove
that the R.\,Thompson group
satisfies the rational Novikov conjecture.

The first (and still very useful) definition of diagram groups (see \cite{Kil},\cite{GuSa97}) was algebraic. From this point of view, every
diagram group $\dg$ is determined by a semigroup presentation
$\pp$ and a distinguished word $w$. One can give an equivalent
topological definition of diagram groups
\cite{GuSa01}. From the topological point of
view, diagram groups are determined by a directed 2-complex $K$
(all edges have directions, every
cell is bounded by two positive paths, the top and the bottom),
and a distinguished positive path $p$.
Diagram groups are
similar to second relative homotopy groups of 2-complexes,
only one needs to consider directed 2-complexes
and homotopies consisting of positive paths only
(we call them {\em directed
homotopies}).

Here is an informal definition of diagram groups (see~\cite{GuSa97}
for details).

Let $\pp=\rav$ be a semigroup presentation where $\Sigma$ is an
alphabet and $\rr$ is the set of defining relations.

Any diagram over $\pp$ is obtained as follows. Start with a
positive (horizontal) path $p$ on the plane labeled by some word
$w$ over $\Sigma$ (that is, a linear oriented labeled graph with
$|w|$ edges which form a path, whose label is $w$). This is a {\em
trivial $(w,w)$-diagram}, and $p$ is the {\em top} and the {\em
bottom} path of this diagram.

Next find a subword in $w$ which is equal to $u$ (or $v$) for some
relation $u=v$ in $\rr$: $p=p'qp''$ where the label of $q$ is $u$
(resp. $v$). Below $p$, draw a path $q'$ labeled by $v$ (resp.
$u$) whose initial and terminal vertices coincide with the initial
and terminal vertices of $q$. The path $q(q')\iv$ must bound a
region on the plane (called a {\em cell}). The result of this
operation is a one-cell diagram whose top path is labeled by $w$
and the bottom path is labeled by the word obtained from $w$ by
replacing $u$ by $v$ (resp. $v$ by $u$). Attaching a new cell to
the bottom path of the diagram, we get a diagram with two cells,
etc. Every diagram $\Delta$ is a plane labeled oriented graph
which tesselates a region of the plane between two positive paths
$\topp(\Delta)$ and $\bott(\Delta)$. If $w$ is the label of
$\topp(\Delta)$ and $w'$ is the label of $\bott(\Delta)$ then
$\Delta$ is called a $(w,w')$-diagram.

Two diagrams are called {\em equal} if there exists an isotopy of
the plane which takes one of the diagrams to the other one.

A diagram is called {\em reduced} if it does not contain dipoles.
A {\em dipole} is a pair of cells such that the bottom path of one
of them coincides with the top path of the other one and these
cells are mirror images of each other. If a diagram contains a
dipole, the two cells forming the dipole can be removed. So every
diagram can be reduced. By the theorem of Kilibarda~\cite{Kil} the
reduced form of every diagram is unique.

Fix a word $w$ and consider the set $\dg$ of all $(w,w)$-diagrams
over $\pp$. One can multiply two diagrams $\Delta_1$ and
$\Delta_2$ in $\dg$ by gluing together $\bott(\Delta_1)$ and
$\topp(\Delta_2)$ and reducing the resulting diagram. This
operation is associative, the trivial $(w,w)$-diagram plays the
role of the identity element, and every diagram $\Delta$ has an
inverse, the mirror image of $\Delta$. Thus $\dg$ is a group which
is called the {\em diagram group over the presentation $\pp$ with
base word $w$}.

Since we are going to use only the algebraic definition of diagram groups, we
do not give here a precise topological definition.
Let us only mention that the directed complex corresponding
to a semigroup presentation is similar to the standard 2-complex
of a group presentation. It has one vertex, one oriented edge for
each generator and one oriented cell for each relation $u=v$ with bottom path
$u$ and the top path $v$. Then
the word $w$ in the algebraic definition of a diagram group turns
into a positive path $w$ in the directed 2-complex, and every
$(w,w)$-diagram is a planar representative of a directed homotopy
from $w$ to $w$. Conversely, every directed $(w,w)$-homotopy is
represented by a $(w,w)$-diagram. The product of homotopies
corresponds to the product of diagrams. Equivalent diagrams
correspond to equivalent (isotopic) homotopies. This allows one to
translate every statement about diagram groups from the algebraic
language to the topological language and back.

The relation between diagram groups and semigroup presentations
(directed complexes) is not rigid. For example, if the
presentation $\pp$ is aspherical, then the diagram groups are
trivial (regardless of the base). On the other hand, presentations
of finite semigroups may correspond to ``large" diagram groups. In
particular, the diagram group corresponding to the presentation
$\la x\mid x^2=x\ra$ of the trivial semigroup is the well known
R.\,Thompson group $F$ (for every base). The directed complex
corresponding to this presentation is the well known {\em Dunce
hat} \cite{Zee} which can be obtained from the triangle

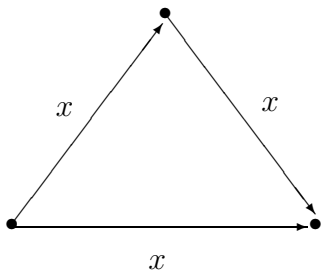
\begin{center}
\unitlength=1mm
\special{em:linewidth 0.4pt}
\linethickness{0.4pt}
\begin{picture}(41.67,34.67)
\put(0.67,5.67){\circle*{1.49}}
\put(41.00,5.67){\circle*{1.49}}
\put(21.00,33.67){\circle*{1.49}}
\put(0.67,5.67){\vector(3,4){20.00}}
\put(21.00,33.67){\vector(3,-4){20.00}}
\put(0.67,5.33){\vector(1,0){39.33}}
\put(7.67,21.00){\makebox(0,0)[cc]{$x$}}
\put(35.00,21.67){\makebox(0,0)[cc]{$x$}}
\put(20.00,0.67){\makebox(0,0)[cc]{$x$}}
\end{picture}
\end{center}

\noindent by gluing all three sides according to their direction.
It is easy to construct other semigroup presentations (directed
complexes) with diagram groups isomorphic to $F$. Nevertheless in
this paper we show that if $F$ appears in a diagram group of a
directed complex (resp. presentation of a semigroup) then Dunce
hat maps into the complex (the semigroup contains an idempotent).
Thus there is a rigid relationship between $F$ and the Dunce hat
(the presentation $\la x\mid x^2=x\ra$).

Recall that the group $F$ can be given by the following infinite
presentation:
$$
\langle\,x_0,x_1,\ldots\mid x_j^{x_i}=x_{j+1}\ (j>i)\,\rangle.
$$
It has also a finite presentation
\be{relF}
\langle\,x_0,x_1\mid x_2^{x_1}=x_3, x_3^{x_1}=x_4\,\rangle,
\ee
where $x_2=x_1^{x_0}$, $x_3=x_2^{x_0}$, $x_4=x_3^{x_0}$ by
definition.

\begin{thm}
\label{th1}
The following conditions are equivalent.
\begin{enumerate}
\item For some word $w$, the diagram group  $\dg$ contains an isomorphic copy of
the R.\,Thomp\-son group $F$.
\item The semigroup given by $\rav$ contains an idempotent.
\end{enumerate}
\end{thm}

A part of this theorem, namely the implication $2\Longrightarrow1$, has been
proved in~\cite[Theorem 25]{GuSa99}. We asked \cite[Problem 2]{GuSa99} whether
the converse is true. Theorem~\ref{th1} gives an affirmative answer to this
question.

The topological formulation of Theorem~\ref{th1} is the following

\begin{thm}
\label{th2}
Let $K$ be a directed complex. Then the following conditions are equivalent.
\begin{enumerate}
\item A diagram group corresponding to $K$ contains an isomorphic copy of
the R.\,Thompson group $F$.
\item The complex $K$ contains a positive non-empty path $t$ which is directly
homotopic to its square.
\item There exists a directed morphism from the Dunce hat to $K$.
\end{enumerate}
\end{thm}

Clearly Theorem \ref{th1} and \ref{th2} are equivalent
We shall prove the theorem in the first formulation.

Recall also that in~\cite[Theorem 24]{GuSa99}, we have proved a
similar rigidity theorem for the restricted wreath product $\zz\wr\zz$.
Similar rigidity theorems might be true for other diagram groups as well.

\section{Proof of the rigidity theorem}

We need one auxiliary geometric fact. Let $\pp$ be a semigroup
presentation and let $\Delta$ be a diagram over $\pp$. For any two
vertices $o'$, $o''$ in $\Delta$ we put $o'\le o''$ whenever there
exists a positive path in $\Delta$ from $o'$ to $o''$. It is easy
to see that the labels of any two positive paths from $o'$ to
$o''$ are equal modulo $\pp$ (see~\cite{GuSa97}). So one can
define the element $\mu(o',o'')$ in the monoid $M$ presented by
$\pp$. This element is represented in $M$ by the label of any
positive path from $o'$ to $o''$.

Recall~\cite{GuSa97} also that for every $(u,v)$-diagram $\Delta$
and $(u',v')$-diagram $\Delta'$ one can define the {\em sum}
$\Delta+\Delta'$ by gluing the terminal vertex of $\topp(\Delta)$
with the initial vertex of $\topp(\Delta')$. The result is a
$(uu',vv')$-diagram. \vspace{2ex}

\begin{lm}
\label{lm1} Let $\pp$ be a semigroup presentation. Let $M$ denote
the monoid presented by $\pp$. Suppose that $\Delta$ is a
$(uv,uv)$-diagram over $\pp$. Let $o_1$ $($resp. $o_2)$ be the
vertex in the top (bottom) path of $\Delta$ that subdivides it
into a product of two paths labeled by $u$ and $v$. Suppose that
$o$ is a vertex in $\Delta$, where $o\le o_1$, $o\le o_2$. If
$\Delta$ is equivalent to a sum of a $(u,u)$-diagram and a
$(v,v)$-diagram, then $\mu(o,o_1)=\mu(o,o_2)$. $($It would be more
precise to write $\mu_\Delta$ but it will always be clear what
diagram we refer to$)$.
\end{lm}
\vspace{1ex}

{\bf Proof.}\ \ Obviously, the reduced form of $\Delta$ is a sum
of a $(u,u)$-diagram and a $(v,v)$-diagram. Suppose that we need
to cancel $m\ge0$ pairs of dipoles in order to reduce $\Delta$. We
prove the claim by induction on $m$. If $m=0$ then the conclusion
is obvious since in this case $o_1=o_2$. Let $m>0$. Cancel a
dipole that consists of two cells $\pi_1$ and $\pi_2$, where the
bottom path of $\pi_1$ coincides with the top path of $\pi_2$. As
a result, we get a diagram $\Delta'$ that can be reduced in $m-1$
step. Let $p_1$ be the top path of $\pi_1$, $p_2$ be the bottom
path of $\pi_2$, and let $p$ be the common boundary of $\pi_1$ and
$\pi_2$.

Suppose first that $o$ is a vertex that does not disappear in
$\Delta'$, that is, $o$ does not belong to $p$ as an inner point.
In this case, for any positive path from $o$ to $o_1$ in $\Delta$,
we can find a positive path from $o$ to $o_1$, which does not
contain $p$ as a subpath (just replace $p$ by $p_1$). The same is
true for positive paths from $o$ to $o_2$. The vertices $o$,
$o_1$, $o_2$ still exist in $\Delta'$, and the elements
$\mu(o,o_1)$, $\mu(o,o_2)$ do not change when we replace $\Delta$
by $\Delta'$. Applying our inductive assumption to $\Delta'$, we
see that these elements are equal there. Thus they are equal for
$\Delta$, too.

Now suppose that $o$ disappears in $\Delta'$. Thus $p$ is
subdivided by $o$ into two paths, say, $q$ and $r$. Let $\bar o$
be the terminal point of $r$. Obviously, any positive path from
$o$ to $o_1$ or $o_2$ begins with $r$. By the previous paragraph,
$\mu(\bar o,o_1)=\mu(\bar o,o_2)$. Now it remains to notice that
$\mu(o,o_i)=\nu\mu(\bar o,o_i)$, where $\nu\in M$ is the element
represented by the label of $r$ ($i=1,2$). This completes the
proof.
\vspace{2ex}

Let $H$ be a group, $y_0,y_1\in H$. Suppose that $y_0$, $y_1$ do
not commute in $H$ and satisfy relations~(\ref{relF}), that is,
$y_1^{y_0y_1}=y_1^{y_0^2}$, $y_1^{y_0^2y_1}=y_1^{y_0^3}$. Since
all proper homomorphic images of $F$ are abelian~\cite{CFP}, it is
clear that $y_0$, $y_1$ generate $F$ as a subgroup of $H$. In this
case, we say that an ordered pair $y_0$, $y_1$ {\em generates $F$
canonically}. We also introduce elements $y_i$ for $i\ge2$ by
$y_i=y_1^{y_0^{i-1}}$.

Let us recall some definitions. We refer to~\cite[Section
15]{GuSa97} for details. Let $\pp=\pres{\Sigma}{\rr}$ be any
semigroup presentation. A $(w,w)$-diagram $\Delta$ over $\pp$ is
called {\em absolutely reduced} provided $\Delta^n$ is reduced for
every $n\ge1$. For any $(w,w)$-diagram $\Delta$ over $\pp$, where
$w\in\Sigma^+$, there exists a word $v\in\Sigma^+$, a
$(w,v)$-diagram $\Psi$ and an absolutely reduced $(v,v)$-diagram
$\bar\Delta$ such that $\Delta=\Psi\bar\Delta\Psi\iv$. One can
decompose $\bar\Delta$ into a sum $A_1+\cdots+A_m$ of spherical
diagrams. Here each nontrivial summand cannot be decomposed into a
sum of spherical diagrams. We also assume that for any $i$ ($1\le
i<m$) at least one of the diagrams $A_i$, $A_{i+1}$ is nontrivial.
The summands $A_i$ ($1\le i\le m$) are called {\em components\/}
of $\bar\Delta$. The number of nontrivial components does not
depend on the choice of $\bar\Delta$. So it can be denoted by
$\comp(\Delta)$.

Let $G=\dg$ be a diagram group. For any $(\bs,v)$-diagram $\Psi$ over $\pp$,
where $v\in\Sigma^+$, we have an isomorphism
$\psi\colon G\to H={\cal D}(\pp,v)$ that takes any diagram $\Delta\in G$ to
$\Psi^{-1}\Delta\Psi$. For any $(\bs,\bs)$-diagram $\Delta$ over $\pp$ we can
construct an isomorphism $\psi$ defined above such that the diagram
$\psi(\Delta)=\Psi^{-1}\Delta\Psi$ will be absolutely reduced. In this case
we will often assume without loss of generality that $\Delta$ is absolutely
reduced up to changing the base of our diagram group.

Suppose that $A=A_1+\cdots+A_m$ and $B=B_1+\cdots+B_n$ are absolutely
reduced diagrams each decomposed into a sum of components. Let $A_i$
($1\le i\le m$) and $B_j$ ($1\le j\le n$) be $(v_i,v_i)$- and
$(w_j,w_j)$-diagrams, respectively. If there exists a $(w,v)$-diagram
$\Gamma$ such that $A=\Gamma^{-1}B\Gamma$, where $v=v_1\ldots v_m$,
$w=w_1\ldots w_n$, then $m=n$ and $\Gamma$ can be decomposed into a
sum $\Gamma_1+\cdots+\Gamma_m$ of $(w_i,v_i)$-diagrams $\Gamma_i$ such
that $A_i=\Gamma_i^{-1}B_i\Gamma_i$ ($1\le i\le m$). Any element $C$ in
the centralizer of $A$ in ${\cal D}(\pp,v)$ can be decomposed into a sum
$C=C_1+\cdots+C_m$, where $C_i$ is a $(v_i,v_i)$-diagram that commutes
with $A_i$ ($1\le i\le m$). If $A_i$ is nontrivial then its centralizer
is cyclic so $A_i$ and $C_i$ belong to the same cyclic subgroup. It is easy
to see that one can change the base in such a way that both diagrams
$A$ and $C$ become cyclically reduced (see~\cite[Theorem 17]{GuSa99}).
\vspace{2ex}

The following theorem is stronger than the implication $1\Longrightarrow2$ in
Theorem~\ref{th1}.

\begin{thm}
\label{th3}
Let $\pp=\rav$ be a semigroup presentation, $w\in\Sigma^+$. If the diagram
group $G=\dg$ contains an isomorphic copy of R.\,Thomp\-son's group $F$, then
the semigroup $S$ presented by $\pp$ contains an idempotent. Moreover, $G$
contains a copy of $F$ if and only if there exist words $w_1,w_2\in\Sigma^*$,
$e\in\Sigma^+$ such that equalities $w=w_1ew_2$, $e^2=e$ hold modulo $\pp$.
\end{thm}

{\bf Proof.}\ \ Suppose that $G=\dg$ contains an isomorphic copy
of $F$. Then there exist $(\bs,\bs)$-diagrams $Y_0$, $Y_1$ over
$\pp$ that generate $F$ canonically.  We assume that the total
number of their components, that is, $\comp(Y_0)+\comp(Y_1)$, is
minimal possible. Note that this number does not change if we
replace $Y_0$, $Y_1$ by their conjugates $Y_0^D$, $Y_1^D$ for any
$(\bs,v)$-diagram $D$ over $\pp$, where $v$ is a nonempty word
over $\Sigma$.

It is easy to see that the element $x_2x_3x_2^{-2}\in F$ commutes
with $x_i$ for all $i\ge3$. So it also commutes with
$x_3x_4x_3^{-2}$. Changing the base $w$, we can assume without
loss of generality that $D_2=Y_2Y_3Y_2^{-2}$ is a cyclically
reduced diagram over $\pp$ decomposed into the sum of components
$A_1+\cdots+A_m$, where $A_i$ is a $(v_i,v_i$)-diagram $(1\le i\le
m)$. Obviously, $D_2$ is nontrivial. (Otherwise $Y_2=Y_3$ and
$Y_0$ commutes with $Y_1$ .) Since $D_3=Y_3Y_4Y_3^{-2}$ is in the
centralizer of $D_2$, we can assume that both diagrams $D_2$,
$D_3$ are absolutely reduced and $D_3=B_1+\cdots+B_m$, where $B_i$
commutes with $A_i$ for all $1\le i\le m$. (Note that the summands
$B_i$ are not necessarily components of $D_3$.)

Suppose that $A_i$ is nontrivial for some $i$. Let it be the $j$th
nontrivial component of $D_2$ counting from left to right. It is
clear that $D_3=D_2^{Y_0}=D_2^{Y_1}$. Thus $D_2$ and $D_3$
conjugate and so they have the same structure of components. Let
$B'$ be the $j$th nontrivial component of $D_3$ counting from left
to right. There are three possible cases: $B'$ is contained in
either 1) $B_i$, or 2) $B_1+\cdots+B_{i-1}$, or 3)
$B_{i+1}+\cdots+B_m$. Clearly, the third case is symmetric to the
second one. So we consider only the first two cases.

{\bf Case 1.}\ \ It is obvious that $B'=B_i$. Hence the
conjugation of $D_2$ by each of $Y_0$, $Y_1$ takes $A_i$ to $B_i$.
This implies that each of the diagrams $Y_0$, $Y_1$ can be
decomposed into a sum of three spherical diagrams with bases
$u_1=v_1\ldots v_{i-1}$, $u_2=v_i$, $u_3=v_{i+1}\ldots v_m$,
respectively. So we have an injective homomorphism $\phi$ from the
Thompson group $F$ (generated by $Y_0$, $Y_1$) to the direct
product ${\cal D}(\pp,u_1)\times{\cal D}(\pp,u_2)\times{\cal
D}(\pp,u_3)$. Denote by $H_k$ the projection of $F$ onto $k$th
factor and let $\psi_k$ be the homomorphism from $F$ onto $H_k$
($k=1,2,3$). The group $F$ embeds into $H_1\times H_2\times H_3$.
Therefore, at least one of the three groups $H_k$ is not abelian.
Then it must be isomorphic to $F$ because all proper homomorphic
images of $F$ are abelian. So let $H_k$ be non-abelian. Let us
show that $k=1$ or $k=3$.

We know that the diagrams $A_i$, $B_i$ belong to the same cyclic
subgroup. By~\cite[Theorem 15.30]{GuSa97}, we may assume that they
belong to the maximal cyclic subgroup $K$ of the diagram group
${\cal D}(\pp,v_i)$. Let us establish that any $(v_i,v_i)$-diagram
$D$ over $\pp$ such that $A_i^D=B_i$, also belongs to $K$. Let $C$
be the generator of $K$. By definition, $A_i$ is nontrivial. So
$B_i$ is also nontrivial and so we have $A_i=C^r$, $B_i=C^s$,
where $r$, $s$ are non-zero integers. We now have $(C^D)^r=C^s$.
So we can apply~\cite[Corollary 15.28]{GuSa97} to conclude that
there is a diagram $C_0$ and some integers $p$, $q$ such that
$C^D=C_0^p$, $C=C_0^q$ and $pr=qs$. Since $C$ generates maximal
cyclic subgroup, we have $|p|=|q|=1$. Thus $C^D=C^{\pm1}$. If
$C^D=C^{-1}$, then $(CD)^2=D^2$. Using the fact that diagram
groups have the unique extraction of roots property (\cite[Section
15]{GuSa97}), we deduce that $C$ is trivial. This is a
contradiction. So $C^D=C$. Hence $D$ belongs to $K$ because  $K$
coincides with its centralizer. Now we can conclude that the
images of $Y_0$, $Y_1$ under $\psi_2$ belong to the same cyclic
subgroup. So $H_2=\psi_2(F)$ is abelian.

We have proved that either $H_1$ or $H_3$ is isomorphic to $F$. It
is obvious that for any diagram $\Delta$ from the subgroup
generated by $Y_0$, $Y_1$, one has
$\sum_{k=1}^3\comp(\psi_k(\Delta))=\comp(\Delta)$. Since
$\psi_2(F)$ is nontrivial, we see that
$\comp(\psi_2(Y_0))+\comp(\psi_2(Y_1))>0$. So for any $k=1,3$ we
have
$\comp(\psi_k(Y_0))+\comp(\psi_k(Y_1))<\comp(Y_0)+\comp(Y_1)$. Now
we can take the value of $k$ such that $\psi_k(F)\cong F$ and
replace the elements of our canonical generating pair $Y_0$, $Y_1$
by their images under $\psi_k$. We get another canonical
generating pair with smaller total number of components. This is a
contradiction, so Case 1 is impossible.
\vspace{1ex}

{\bf Case 2.}\ \ Let $B'$ be contained in $B_1+\cdots+B_{i-1}$ as a
subdiagram. We have $B_1+\cdots+B_{i-1}=\Xi_1+B'+\Xi_2$ for some
spherical diagrams $\Xi_1$, $\Xi_2$. Let $z$ be the base of the diagram
$B_{i+1}+\cdots+B_m$ and let $t$ be the base of $\Xi_2+B_i$. Obviously,
$t$ is nonempty because it has a terminal segment $v_i$. We will show
that $t^2=t^3$ modulo $\pp$ so $e=t^2$ represents an idempotent in $S$.
It will be also clear that $w$ belongs to the two-sided ideal in $M$
generated by $e$, where $M=S^1$ is the monoid presented by $\pp$.

Let $D$ be $Y_0$ or $Y_1$. We use the fact that $D_2^D=D_3$. Each
of the diagrams $D_2$, $D_3$ is a sum of components. According to
the above description, $D$ can be naturally decomposed into a sum
of $m$ diagrams (not necessarily spherical) such that the
conjugation by the $k$th summand $(1\le k\le m$) takes $A_k$ (the
$k$th component of $D_2$) to the $k$th component of $D_3$ (recall
that this component may not coincide with $B_k$). Then $A_i$, the
$j$th nontrivial component of $D_2$, is taken to $B'$, the $j$th
nontrivial component of $D_3$. The bases of diagrams to the right
of $A_i$, $B'$ in $D_2$ and $D_3$, respectively, are $z$ and $tz$.
This means that $D$ is a sum of an $(xt,x)$-diagram and a
$(z,tz)$-diagram, where $x$ is the base of $\Xi_1$.

Note that $x_2x_3x_2^{-2}$ commutes with $x_3$. So $Y_3$ belongs
to the centralizer of $D_2$. Hence $Y_3$ is a sum of an
$(xt,xt)$-diagram and a $(z,z)$-diagram. The diagram $$
\Delta\equiv Y_0^{-1}\circ Y_0^{-1}\circ Y_1\circ Y_0\circ Y_0, $$
equivalent to $Y_3$, has the following structure:

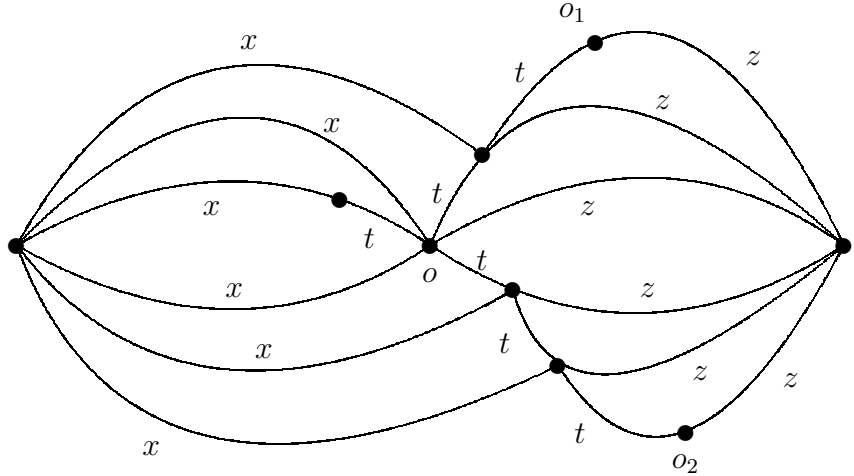
\begin{center}
\unitlength=1mm
\special{em:linewidth 0.4pt}
\linethickness{0.4pt}
\begin{picture}(112.00,81.00)
\put(1.00,31.00){\circle*{2.00}}
\put(56.00,31.00){\circle*{2.00}}
\put(111.00,31.00){\circle*{2.00}}
\bezier{260}(1.00,31.00)(31.00,48.00)(56.00,31.00)
\bezier{264}(56.00,31.00)(86.00,49.00)(111.00,31.00)
\bezier{260}(1.00,31.00)(31.00,14.00)(56.00,31.00)
\bezier{264}(56.00,31.00)(83.00,13.00)(111.00,31.00)
\put(44.00,37.00){\circle*{2.00}}
\put(67.00,25.00){\circle*{2.00}}
\bezier{348}(1.00,31.00)(24.00,1.00)(67.00,25.00)
\bezier{296}(67.00,25.00)(74.00,0.00)(111.00,31.00)
\put(73.00,15.00){\circle*{2.00}}
\bezier{432}(1.00,31.00)(17.00,-12.00)(73.00,15.00)
\bezier{304}(73.00,15.00)(89.00,-10.00)(111.00,31.00)
\put(90.00,6.00){\circle*{2.00}}
\bezier{352}(1.00,31.00)(34.00,65.00)(56.00,31.00)
\bezier{376}(56.00,31.00)(71.00,68.00)(111.00,31.00)
\put(63.00,43.00){\circle*{2.00}}
\bezier{384}(1.00,31.00)(24.00,72.00)(63.00,43.00)
\bezier{400}(63.00,43.00)(87.00,81.00)(111.00,31.00)
\put(78.00,58.00){\circle*{2.00}}
\put(32.00,58.00){\makebox(0,0)[cc]{$x$}}
\put(68.00,54.00){\makebox(0,0)[cc]{$t$}}
\put(99.00,56.00){\makebox(0,0)[cc]{$z$}}
\put(43.00,47.00){\makebox(0,0)[cc]{$x$}}
\put(27.00,36.00){\makebox(0,0)[cc]{$x$}}
\put(48.00,32.00){\makebox(0,0)[cc]{$t$}}
\put(30.00,25.00){\makebox(0,0)[cc]{$x$}}
\put(57.00,38.00){\makebox(0,0)[cc]{$t$}}
\put(87.00,50.00){\makebox(0,0)[cc]{$z$}}
\put(77.00,36.00){\makebox(0,0)[cc]{$z$}}
\put(85.00,25.00){\makebox(0,0)[cc]{$z$}}
\put(104.00,13.00){\makebox(0,0)[cc]{$z$}}
\put(92.00,14.00){\makebox(0,0)[cc]{$z$}}
\put(19.00,4.00){\makebox(0,0)[cc]{$x$}}
\put(76.00,6.00){\makebox(0,0)[cc]{$t$}}
\put(34.00,17.00){\makebox(0,0)[cc]{$x$}}
\put(66.00,18.00){\makebox(0,0)[cc]{$t$}}
\put(63.00,29.00){\makebox(0,0)[cc]{$t$}}
\put(56.00,27.00){\makebox(0,0)[cc]{$o$}}
\put(75.00,62.00){\makebox(0,0)[cc]{$o_1$}}
\put(90.00,2.00){\makebox(0,0)[cc]{$o_2$}}
\end{picture}
\end{center}

\noindent Here $o$ is the vertex in $Y_1$ that subdivides it into
the sum of an $(xt,x)$- and a $(z,zt)$-diagrams. By $o_1$ ($o_2$)
we denote the vertex on the top (bottom) path of $\Delta$ that
subdivides this path into a product of paths with labels $xt$ and
$z$. Clearly, there is a path in $\Delta$ from $o$ to $o_1$
labeled by $t^2$ and there is a path in $\Delta$ from $o$ to $o_2$
labeled by $t^3$. Applying Lemma \ref{lm1}, we conclude that
$t^2=t^3$ modulo $\pp$. (It is obvious that $w$ belongs to $Mt^2M$
as an element in $S$.) \vspace{1ex}

The converse is proved in~\cite[Theorem 25]{GuSa99}.
\vspace{1ex}

The proof is complete.

\begin{rk}
{\rm
Given a finite semigroup presentation $\pp$ and a word $\bs\in\Sigma^+$,
we cannot decide algorithmically whether the diagram group $\dg$ contains
$F$ as a subgroup. Indeed, the property of a finitely presented semigroup
not to have an idempotent, is a Markov property. Let $a$, $b$ be new
letters that do not belong to $\Sigma$. Adding them to $\Sigma$ and
adding relations of the form $ax=a$, $xb=b$ ($x\in\Sigma$), we get a
new semigroup presentation $\qq$. The diagram group ${\cal D}(\qq,ab)$
contains $F$ as a subgroup if and only if $S$ has an idempotent, where
$S$ is the semigroup presented by $\pp$. This is clear because all
idempotents in the semigroup presented by $\qq$ are represented by
words over $\Sigma$ and $ab$ belongs to the two-sided ideal generated by
any word over $\Sigma$.
}
\end{rk}

\noindent Victor Guba\\ Vologda State Pedagogical University,\\
S.~Orlov Street 6,\\ Vologda, 160600\\ Russia\\
guba{@}uni-vologda.ac.ru\\

\medskip

\noindent Mark Sapir\\ Vanderbilt University, \\Nashville, TN
37240, U.S.A.\\ msapir{@}math.vanderbilt.edu

\end{document}